\documentstyle[12pt]{article}
\makeatletter\@addtoreset{equation}{section}\makeatother

\def\Mat{{\rm Mat}}
\def\d{\delta}
\def\e{\epsilon}

\def\dim{{\rm dim }}
\def\Tr{{\rm Tr}}

\def\K{{\bf K  }}

\def\ind{{\rm ind }}
\def\deg.tr.{{\rm deg.tr. }}
\def\rank{{\rm rank }}
\title{ On a new class of the commutative subalgebras 
of the maximal Gel'fand-Kirillov dimension in the 
universal enveloping algebra of a simple Lie algebra }
\author {T. Skrypnyk \footnote{The research described in this publication was possible in part by the Award 
number UP1-309 of the U.S. Civilian Research and Development Foundation 
(CRDF) for independent states of the former Soviet Union.}}

\date{ }
\begin{document}
\thispagestyle{empty}
\maketitle
 {\it Bogolyubov Institute for Theoretical Physics,
Metrologichna st.,14 b., 252143, Kiev, Ukraine.
e-mail : tskrypnyk@elc3.imath.kiev.ua }
\begin {abstract}
New commutative subalgebras of the maximal Gel'fand-Kirillov dimension
in the universal enveloping algebras  of the classical Lie algebras $gl(n)$ and $so(n)$ are constructed. In the case of $sp(n)$
Gel'fand-Tsetlin algebra is extended to a maximally commutative one. 
\end{abstract}
\newpage
\section{Introduction}
 Commutative subalgebras of maximal degree of transcendency (it is known \cite{GK} that for commutative associative algebra the degree of transcendency coinsides with Gel'fand- Kirillov dimension) in the universal enveloping algebras are very important both
from the point of view of the representation theory and the theory of quantization of completely
integrable systems. In the framework of the latter these subalgebras play the role of the algebras of quantum integrals of motion.
In both theories a possibility to diagonalize the generators of  these algebras is very important.
It may yield the basis in the space of representation which consists of common eigenvectors of all generators of these algebras. Examples of such the construction are known: so-called Gel'fand-Tsetlin algebras  consisting of the Casimir elements of a chain of embedded  subalgebras.
Bases of their common eigenvectors are celebrated Gel'fand-Tsetlin bases  \cite{GZ}-\cite{G}.
Except for the special case of \cite{Ya} there were no commutative subalgebras of the maximal degree of transcendency other than Gel'fand-Tsetlin algebras known in the universal enveloping algebras of the simple Lie algebras.

 In this paper we find new commutative subalgebras in the  universal enveloping algebras $U(g)$ of the simple Lie  algebras $g$  that have maximal degree of transcendency and include any given commutative subalgebra $C_{g_A}\subset U(g_A)$  for the next  pairs $(g,g_A)$ of  Lie algebras and their subalgebras: $( gl(n), gl(n-2)+gl(2) ), (so(n), so(n-2)+so(2))$ and $ (sp(n), sp(n-1)+gl(1))$.
 Idea that lead to the discovery of such the algebras comes out from the theory of integrable dynamical systems. As it is known \cite{MF}- \cite{MF2} there is sufficiently large (multiparametric) class of commutative with respect to the Lie-Poisson brackets  subalgebras in the algebra of the polinomial functions $P(g^*)$ on the dual to the Lie algebra space $g^*$.
Algebra of polinomial functions on $g^*$ can be viewed as a graded algebra of the filtered algebra $U(g)$ and the Lie-Poisson brackets as the graded part of the commutator in $U(g)$.
Equality to zero of the Lie-Poisson brackets of the elements of $P(g^*)$ is equivalent to the equality to zero of the higher order in the elements of $g$ part of the commutator (but not of the whole commutator!) of the correspondent elements of $U(g)$.
This means that, in general, commutative with respect to the Lie-Poisson brackets subalgebra in $P(g^*)$ may not have commutative counterpart in $U(g)$.
But if  a commutative subalgebra in $P(g^*)$ is an algebra of integrals of some dynamical system of the Euler-Arnold type then from the physical consideration it follows, that corresponding quantum dynamical system should have the same algebra of the integrals, meaning  that counterpart of commutative algebra in $P(g^*)$  have also to be commutative subalgebra in $U(g)$.
For the case when underlying dynamics  is the dynamics of the generalized rigid body the corresponding algebra of the integrals coincides with the Mishchenko-Fomenko algebra \cite{MF}-\cite{MF2}. Although it was conjectured by Vinberg \cite{V} that the Mishchenko-Fomenko algebra could be "lifted" to the universal enveloping algebra without the lost of the commutativity,  he managed to prove this fact only for its special subalgebra, characterized as having a second order in the elements of the Lie algebra $g$.
We prove (Theorems 1-2) the commutativity in $U(g)$ of the other analogue of the subalgebra of the Mishchenko-Fomenko algebra. Contrary to the case considered by Vinberg, our subalgebra in the $U(g)$  is not "homogeneous"  in the elements of $g$ but is linear in the "parameters of the shift" $(A)_{i,j}$, i.e. additional parameters which are important ingradients in the definition of the Mishchenko-Fomenko algebras \cite{MF}-\cite{MF2} \footnote{For definition of the Mishchenko-Fomenko algebras see also section 6.}.
We use it  to construct commutative subalgebra in $U(g)$,
including its  generators,   the Casimir elements of $U(g)$ and  generators of  a maximally commutative subalgebra in $U(g_A)$, where  $g_A$ is a stabilizer of  a fixed matrix $A$ in $g$, in the set of generators of the commutative subalgebra  in $U(g)$. Using results of \cite{Br}, \cite{BS}, we prove (Theorem 3) that in the cases when  "shift" matrix $A$ is of the second  rank these commutative subalgebras  in $U(g)$ have the maximal degree of transcendency equal to $(\dim g +\ind g)/2$.

Importance of these commutative subalgebras is explained  
by their  maximal degree of transcendency and the fact that
they contain the maximal number of generators that commute with the described above Lie subalgebras $g_A\in g$.
These both facts together mean that they should separate all vectors of the isomorphic (multiple) representations of the Lie subalgebra $g_A$ that emerge after a restriction of a representation of the Lie algebra $g$ to its subalgebra $g_A$. \footnote {See  \cite{Ya} for a "physical" example of the representations of  $su(2,2)\subset gl(4,{\bf C})$.}
 
   The structure of our article is following: in section 2 we prove commutativity of the described above subalgebras in $U(g)$, in section 3 we prove the maximallity of the degree of transcendency for the constructed commutative subalgebras, and   in the last sections  we give the proofs of the main theorems.

\section {Commutativity}
 Here  we will present some commutative subalgebras in the universal enveloping algebra,
that will be used in the next section to construct maximal commutative subalgebras.
 All considered algebras in this section  are over a field  $\K$  of characteristic zero.

{\bf 2.1} Let us consider the case of $gl(n)$ algebra first.
Let $X_{i,j}, i,j\in I$ where
$I=(1,2,...,n)$, be the basis in this algebra
with the commutation relations:

$[X_{i,j},X_{k,l}]=\d_{k,j} X_{i,l} -\d_{i,l} X_{k,j}$.

Universal enveloping algebra $U(gl(n))$ consists of formal polinomials in the elements $X_{i,j}$.
Let us define the following elements of  $U(gl(n))$:

 $ (X^M)_{i,j}=\sum\limits_ { i_1,...i_{M-1} \in I} X_{i,i_1}X_{i_1,i_2}...X_{i_{M-1},j}$.

It is known \cite{Ba}, that the elements:
$(X^M)=\sum\limits_{i\in I}(X^M)_{i,i}$, $M\in (1,2,...,n)$
generate the center of the universal enveloping algebra.
For arbitrary numerical  matrix $A\in{\Mat (n,\K)}$ and positive integer $M$  we can introduce the following elements of $U(g)$:

$(AX^M)=\sum\limits_{i,j\in I} A_{j,i}(X^M)_{i,j}$.

They  are "quantum " analogues of the generators of the subalgebra of  the Mishchenko-Fomenko algebra linear in 
the "shift" parameters $A_{j,i}$.

Let $g_A\subset gl(n)$ be the subspace in $g$ spanned over  $\K$ by the elements
   $ (BX)=\sum\limits_{i,j\in I}B_{j,i}X_{i,j}$, 
where numerical matrix $B$ belongs to the stabilizer of $A$ in ${ \Mat (n,\K)}$, i. e. $[B,A]=0$.
It is easy to show that 
$[(BX), (CX)]=([B,C]X)$,
i. e. that $g_A$ is a subalgebra in $g$ isomorphic to the stabilizer of matrix $A$ in $\Mat (n,\K)$.

Now we are able to formulate the theorem.

{\large Theorem 1.}
{\it Let $F_{A}$ be a subalgebra in $U(gl(n))$ generated by $(X^N)$ and $(AX^M)$ , $M,N\in{\bf{ Z_{+}}}$.
Then

(i)  $F_A$ is a commutative subalgebra in $U(gl(n))$.

(ii)  $g_A$  centralizes $F_A$ in $gl(n)$.}

{\it Idea of proof}.
Proof of the statement (ii) of  Theorem 1 is based on the fact that elements $(X^M)_{i,j}$ are tensorial \cite{Ba} . We prove statement (i) using the fact that higher order in the elements of $gl(n)$ part of the commutator $[(AX^M),(AX^N)]$  is trivial due to the commutativity of  $(AX^M)$ and $(AX^N)$ in $P(gl(n)^*)$, the observation that the other part of the commutator could be written in the form of the sum of the communators of  $(AX^K)$ and $(AX^L)$ where $K<M$ and the obvious fact that $(AX)\in g_A$.
  See section 4 for detailed proof.

{\bf 2.2} In this subsection we  consider the case of the algebras $so(n)$ and $sp(n)$  which will be denoted by $g$.
It will be convenient to treat them uniformly, so that  the algebras $so(2n)$ and $sp(n)$ be considered as subalgebras of $gl(2n)$ and algebra $so(2n+1)$ as a subalgebra of $gl(2n+1)$ respectively.
Let $X_{i,j}, i,j\in I$ be the basis in these algebras, $I=(-n, ...,-1,0,1,...,n)$ for the case of $so(2n+1)$, $I=(-n,...,-1,1,...,n)$ for the case of $sp(n)$ and $so(2n)$, with the commutation relations:

$[X_{i,j},X_{k,l}]=\d_{k,j} X_{i,l} -\d_{i,l} X_{k,j}+\e_i \e_j(\d_{j,-l}X_{k,-i}-\d_{k,-i}X_{-j,l })$

and additional property 
$X_{i,j}= - \e_i \e_j X_{-j,-i }$
(here $\e_j=1$ for the case $g=so(n)$ and $\e_j=$sgn $j$ for the case $g=sp(n)$).
Universal enveloping algebra $U(g)$ consists of formal polinomials in the elements $X_{i,j}$.
Let $ (X^M)_{i,j}$, $(X^M)$, $(AX^M)$ be  elements of $U(g)$  and $F_A\subset U(g)$, $g_A\subset g$ be subalgebras  of $U(g)$ defined like in previous subsection. 
Then the following theorem analogous to  Theorem 1 holds:

{\large Theorem 2.}
{\it If  matrix A satisfies condition
$A_{i,j}=\pm \e_i \e_j A_{-j,-i} $
then

(i)  $F_A$ is a commutative subalgebra in $U(g)$.

(ii)  $g_A$ coincides with the centralizer of $F_A$ in $g$.}

{\it Idea of proof}.
Proof of Theorem 2 is based on the same idea as the proof of  Theorem 1 (see the end of  the previous subsection).
See section 5  for the detailed proof.

{ \it Remark 1.}
The type of the  realization of algebras $so(n)$ and $sp(n)$ used in this section is not of any importance. For example we could use standard realization of the algebra $so(n)$
with  the basis
$X_{i,j} ,i,j\in (1,2,...,n)$, $X_{i,j}= -X_{j,i }$
that satisfies  the commutation relations:
$[X_{i,j},X_{k,l}]=\d_{k,j} X_{i,l} -\d_{i,l} X_{k,j}+\d_{j,l}X_{k,i }-\d_{k,i}X_{j,l }$.

Proof of  Theorem 2 in this standard  realization is the same as in the previous case.
Important in the above construction is the treatment of the $so(n) (sp(n))$ algebras as subalgebras of $gl(n) (gl(2n))$ that are  stable under the action  of the involutive automorphism $\sigma$.
For elements $(AX^N)\in U(gl(n))$  to remain  commutative after the restriction to the  universal enveloping algebra of these subalgebras for arbitrary realization of the automorphism $\sigma$, we have to require equality
$\sigma (A)=\pm A$,
which in our case has the form:
$A_{i,j}=\pm \e_i \e_j A_{-j,-i} $.

{\it Remark 2.}
  If $A_{i,j}= + \e_i \e_j A_{-j,-i} $
( $\sigma (A)= -A$ in general)
 then "quantum integrals" obtained above correspond to the classical integrals of "normal series" \cite {MF2}. In the next sections we will consider only subalgebras constructed with the help of the matrix $A$, that satisfies the condition
 $A_{i,j}= - \e_i \e_j A_{-j,-i} $, meaning that numerical matrix $A$ belong to the same subalgebra of  $gl(n)$ as the considered Lie algebra.

\section {Maximality}

{\bf 3.1} Algebras $F_A$ constructed in the previous section for all
 classical Lie algebras $g$ \ obviously have
 degrees of transcendency not higher than 2 \rank $g$ -1. Nevertheless statement (ii)
 of  Theorem 1 and Theorem 2 enables us to extend them to the commutative algebras having maximal degree of transcendency.
 Indeed, according to this statement we can use in the construction of maximal commutative
 subalgebra in $U(g)$ not only generators of $F_A$ but also some commutative
 subalgebra in $U(g_A)$. In the case of sufficiently singular matrix $A$
  we may hope to find  in $U(g_A)$ "large " commutative subalgebra,
 that will "complete" $F_A$ to the one having maximal degree of transcendency.
 Indeed, the following theorem is true.

 {\large Theorem 3. }
 {\it Let $g$ be simple matrix Lie algebra over the field of real or complex numbers, $A\in g$ be a constant semisimple matrix of the second matrix rank, $g_A$  and $F_A$ be as in section 2.
 Let $C_{g_A}$  be some commutative subalgebra in $U(g_A)$,
$F (g,C_{g_A})$ be a commutative subalgebra in $U(g)$ generated by $F_A$ and $C_{g_A}$.

If
 $\deg.tr.$ $ C_{g_A}=(\dim g_A+\ind g_A)/2$
then
$\deg.tr. $ $F (g,C_{g_A}) =(\dim g +\ind g)/2$.}
\footnote{The abbreviation $deg.tr.$  means "the degree of transcendency ".}

{\it Idea of proof}. To prove this Theorem we use commutativity of the algebras $F_A$ and $C_{g_A}$, provided by  Theorems 1-2 and the assumption of  Theorem 3, and  the results of Vinberg \cite{V} which enables us to consider  problem of finding the degree of transcendency of our subalgebra in $P(g^*)$ instead $U(g)$. Hence the problem of finding the degree of transcendency is reduced to the problem of finding of the number of functionally independent  generators of  the algebra $F (g,C_{g_A})$ in $P(g^*)$. That permits us to use results of \cite{Br} about a number of functionally independent generators of the Mishchenko-Fomenko algebra in the case of a singular covector of the shift and results of \cite{BS} concerning their explicit form.
See section 6 for detailed proof.

{\it Remark 1.}
The condition that numerical matrix $A$ is of matrix rank 2  is essential. Indeed, if its matrix rank exceeds 2 then from the dimensional considerations follows that $\deg.tr. F (g,C_{g_A}) < (\dim g +\ind g)/2$ for arbitrary  commutative algebra $C_{g_A}\in U(g_A)$ \cite{BS}. On the contrary, it is possible to show that in the case of matrix $A$ of rank 1 elements $(AX^N)$ are functionally dependent on  the Casimir elements of $U(g)$ and $U(g_A)$. So in the last case they will not contribute to the degree of transcendency of the constructed commutative subalgebras.

This theorem enables us to construct the commutative algebras of maximal degree of transcendency  in the  universal enveloping  algebra of an arbitrary simple Lie algebra. Indeed, if there is no known  commutative subalgebra in $U(g_A)$ with the maximal Gel'fand-Kirillov dimension, as in the case of $g=sp(n)$, $g_A=sp(n-1)+gl(1)$, we can use the same method of constructing  commutative subalgebras in the  universal enveloping algebra for each simple component of the $g_A$. Repeating this procedure we at last will come to the situation when all components of  stabilizer satisfy conditions of Theorem 3 (for example when they have $\rank $ one as a simple Lie algebras).

In the general situation in order to obtain "large" commutative subalgebras in the universal enveloping algebra one can combine "method of the  chain of subalgebras" with the method of " the shift of the argument", or to be more precise, its special case connected with special singular "constant covector $A$", varying its degree of singularity from 1 (except the case of symplectic algebra) to 2 for each successive pair of an algebra and its subalgebra in the chain.
 In the next subsections we will consider in the detail  this construction for all  classical simple (reductive) Lie algebras over the field ${\bf C}$ of a complex numbers.
We will use the following easily proved

{\large Lemma 1. }
{\it Let $A\in g\subset Mat(n,\K)$ be numerical semisimple matrix of the second matrix rank and the basic field $\K$ be the field of the complex numbers. 
If $g$ is one of the classical simple Lie algebras of the type $A_n,
B_n, C_n, D_n$ then $g_A$ will contain simple Lie subalgebra of the type  $A_{n-2},B_{n-1}, C_{n-1}, D_{n-1}$ correspondently.}

{\bf 3.2} Let us start with the simplest case $g= gl(n)$ ( or, equivalently, $g=sl(n)$).
 The next corollary of  Theorem 3 and Lemma 2 holds.

 {\large Corollary 1.}
 {\it Let us consider the chain of subalgebras in $gl(n)$ of the type

$gl(n)\supset gl(n-k_1)\supset gl((n-k_1)-k_2))\supset ...\supset gl(1)$.

If $k_i=1,2$   for any index $i$, then one can associate with  this chain a commutative subalgebra of the maximal degree of transcendency (equal to $n(n+1)/2$) in $U(gl(n))$ that consists of the Casimir elements of  the each subalgebra in the chain, and  if subalgebra next to $gl(i)$ in the chain is $gl(i-2)$, elements of $U(gl(n))$ of the type $(A_{(i)}X_{(i)}^N), N<n$.
  Index $(i) $ here means that $X_{(i)}$ belongs to corresponding subalgebra $gl(i)$, $A_{(i)}\in Mat(i)$  is  complex semisimple matrix of the second matrix rank.}

{\it Example.}
 Let $g=gl(n)$; semisimple rank two  matrix $A$ be of the type
$A=diag(a_1,a_2,0,...,0)$, where $ a_i \in {\bf C}, a_i\not= o.$
 Then $g_A=gl(n-2)+gl(1)+gl(1)$, if $a_1\not= a_2$
or $g_A=gl(n-2)+gl(2)$ if $ a_1 = a_2 .$
In both cases we can construct commutative subalgebras in  $U(gl(n) )$  with degree of transcendency equal to $n(n+1)/2$.
For its independent generators we can take:

1)  $n=\rank$ $gl(n)$ independent Casimir elements of type $(X^k)$, $k\in( 1,2,...,n)$.

2)  $n-1$ independent elements  of $U(gl(n))$ of the type $(AX^m)$, $m\in (1,2,...n-1)$.

3) $ (n-1)(n-2)/2$ generators of the Gel'fand-Tsetlin algebra of $gl(n-2)$.

In the second case only we can use instead of the Gel'fand-Tsetlin algebra of $gl(n-2)$ the Gel'fand-Tsetlin algebra of $gl(n-2)+gl(2)$. In this case we have to exclude from the generators of the whole commutative algebra generators $(AX), (AX^2)$ which would be no more algebraically independent. In the special case $n=4$ such the commutative algebra was considered by Yao \cite{Ya} \footnote{In the case  considered by Yao additionally to the Casimir elements of $gl(4)$ and Gel'fand-Tsetlin algebra of $gl(2)+gl(2)$ only one $ Gl(2)\times Gl(2)$ invariant element of $U(gl(4))$ was needed. }.

 {\it Remark 1.}
If $k_1=2$ and complex semisimple matrix $A=A_{(n)}$ of rank two has two equal non zero eigenvalues then all commutative subalgebras, constructed in the Corollary 1, centralize  Lie subalgebra $g_A=gl(n)+gl(2) $.

{\it Remark  2. }
It is interesting to notice that for all commutative subalgebras constructed in the Corollary 1, number of the generators, that belong to the given filtered component of the  universal enveloping algebra, is fixed and is the same as in the Gel'fand-Tsetlin algebra.

 {\bf 3.3}  The case $g=so(n)$ repeats almost completely the case of the $gl(n)$.
 The following corollary of  Theorem 3 and Lemma 1  holds.
 
 {\large Corollary 2.}
 {\it Let us consider the chain of subalgebras in $so(n)$ of the type

$so(n)\supset so(n-k_1)\supset so((n-k_1)-k_2))\supset ...\supset so(2)$.

If $k_i=1,2$ for any index $i$, then one can associate with this chain the commutative subalgebra in $U(so(n))$ of the maximal degree of transcendency equal to $((n-1)n/2+[n/2])/2$, that consists of the  Casimir
elements of  the each subalgebra in the chain  and,  if a subalgebra next to $so(i)$ in the chain is $so(i-2)$, elements of $U(so(n))$ of the type $(A_{(i)}X_{(i)}^N)$, $N<n$ ,$N$ is odd.
Here index $(i)$  means that $X_{(i)}$ belongs to the corresponding subalgebra $so(i)$ and $A_{(i)}\in so(i) $  is  complex semisimple matrix of the second matrix rank. }

 {\it Remark 1.}
If $k_1=2$ then all the commutative subalgebras constructed above centralize
$so(n-2)+so(2)$.

 {\it Remark 2.}
Contrary to the case the $g=gl(n)$ commutative subalgebras  corresponding to the different chains of the subalgebras have different number of the generators, belonging to the given filtered component of the  universal enveloping  algebra. In particular, if we take $k_i =2$, for all indices $i$, we'll  obtain commutative subalgebra in the  universal enveloping  algebra that contains $[n/2]$ linear in elements of $so(n)$ generators. To  make comparison we remind that Gel'fand-Tsetlin algebra for this case has only one generator linear in the elements of the $so(n)$.

{\bf 3.4}  The case $ g=sp(n)$ differs from the considered above cases of $gl(n)$ and $so(n)$. In this case neither for $sp(n)$ nor for any of its subalgebra $sp(n-k)$, $k\in( 1,2,...,n-2)$ Gel'fand -Tsetlin algebra is of maximal degree of transcendency equal to $(n-k)(n-k+1)$. That is why our construction gives only one commutative subalgebra in the $U(sp(n))$ of the maximal degree of transcendency. This commutative subalgebra will be the extended Gel'fand -Tsetlin algebra . This is expressed in the following Corollary:

 {\large Corollary 3.}
{\it Let us consider the following chain of subalgebras in $sp(n)$

$sp(n)\supset sp(n- 1)\supset sp(n- 2)\supset  ...\supset sp(1)  \supset gl(1)$.

A commutative subalgebra in $U(sp(n))$ of the maximal degree of transcendency
equal to $n(n+1)$, that  consists of the Casimir elements of  the each subalgebra in the chain and elements of $U(sp(n))$ of the type $(A_{(i)}X_{(i)}^k)$, $k$ is odd, $k<2i$, is associated with  this chain. Index $(i)$  means that  $X_{(i)}$ belongs to correspondent subalgebra $sp(i)$and $A_{(i)}\in sp(i)$  is complex semisimple matrix of the second rank.}

\section {Proof of  Theorem 1.}
It is easy to verify \cite {Ba} that elements  $ (X^M)_{i,j }$ are "tensorial" :

$\ [X_{i,j}, (X^M)_{k,l}]=\d_{k,j} (X^M)_{i,l} -\d_{i,l} (X^M)_{k,j}.$

This yields the validity of the statement (ii).
Indeed, from the tensorial nature of the elements $(X^N)_{i,j}$  follows that
$[(BX),(AX^N)]=([A,B]X^N)$
Hence $(BX)$ and $(AX^N)$ commute if and only if $[A,B]=0$, i. e. if $(BX)\in g_A$. This proves (ii).

Now we will prove part (i) of the theorem.
It is evident that Casimir elements  of the type $(X^K)$ commute among themselves and with all elements of $U(g)$.
So to prove the statement it is necessary to show that elements of the type $(AX^M)$ and $(AX^N)$ commute.
We will use the following Proposition.

 {\large Proposition 1.} 
$[(X^M)_{i,j},(X^N)_{k,l}]=\sum\limits_{S=1}^{M}((X^{M+N-S})_ {i,l} ( X^{S-1}) _{k,j} - (X^{S-1})_{i,l} (X^{N+M-S})_{k,j}).$

{\it Proof .}
Using the Leibnitz rule  we obtain

 $[ (X^M)_{i,j}, (X^N)_{k,l}]=\sum\limits_{P=1}^{M}\sum\limits_{u,v=1}^{n}((X^{P-1})_ {i,u}[X_{u,v},(X^N)_{k,l}] (X^{ M-P})_{v,j})$.

Using tensorial nature of the elements $(X^N)_{k,l}$  we will obtain the statement of the proposition.

Now we are able to prove (i).
We will use inductive proof.
Let us prove that from
$[(AX^{K}), (AX^N)]=0,$ for all $N\in{\bf  Z_+}, K<M-1 $
follows that
$[(AX^M), (AX^N)]=0,$ for all $N\in {\bf Z_+} $.
This inductive assumption will be proved by the following Proposition:

{\large Proposition 2.} 
$[(AX^M), (AX^N)]= \sum\limits_{S=1}^{M} \sum\limits_{P=1}^{S-1} [(AX^{P-1}),(AX^{M+N-P-1})]$.

{\it Proof .} 
Transforming the left-hand-side of the equality we obtain

$\sum\limits_{i,j,k,l\in I} A_{j,i}A_{l,k} [(X^M)_{i,j},(X^N)_{k,l}]
=\sum\limits_{i,j,k,l\in I} A_{j,i}A_{l,k} \sum\limits_{S=1}^{M}(X^{M+N-S})_ {i,l}  (X^{S-1}) _{k,j} - (X^{S-1})_{i,l} (X^{N+M-S})_{k,j}$

$=\sum\limits_{i,j,k,l\in I} A_{j,i}A_{l,k}(\sum\limits_{S=1}^{M}(X^{M+N-S})_ {i,l}  (X^{S-1}) _{k,j} - (X^{N+M-S})_{k,j} (X^{S-1})_{i,l} - [(X^{S-1})_{i,l}, (X^{N+M-S})_{k,j}])$

$=\sum\limits_{S=1}^{M}((AX^{S-1}AX^{M+N-S})- (AX^{S-1}AX^{M+N-S})- \sum\limits_{i,j,k,l\in I} A_{j,i}A_{l,k}[(X^{S-1})_{i,l},(X^{N+M-S})_{k,j}])$

$= - \sum\limits_{S=1}^{M}\sum\limits_{i,j,k,l\in I} A_{j,i}A_{l,k}[(X^{S-1})_{i,l},(X^{N+M-S})_{k,j}] =-\sum\limits_{i,j,k,l\in I} A_{j,i}A_{l,k} (\sum\limits_{S=1}^{M} \sum\limits_{P=1}^{S-1}(X^{M+N-P-1})_ {i,j}\cdot$

$\cdot ( X^{P-1}) _{k,l}- (X^{P-1})_{i,j} (X^{N+M-P-1})_{k,l})=
\sum\limits_{S=1}^{M}\sum\limits_{P=1}^{S-1} [(AX^{P-1}),(AX^{M+N-P-1})]$.

Now to prove the theorem it necessary only to show, that
$[(AX^M), (AX^N)]=0,$ for all $N\in{\bf Z_+}$, and for $M=0,1$.
 For $M=0$ it is evident, because in this case $(AX^M)=\Tr A$.
Validity of this when $M=1$ follows from the fact, that  $(AX)\in g_A$.
That proves the  theorem.

\section{Proof of Theorem 2}

Part (ii) of  Theorem 2 is proved like this of  Theorem 1 (see the previous section).
Let us prove part (i). 
Like in the case of the previous theorem to prove the statement it is necessary only to show that elements of the type $(AX^M)$ and $(AX^N)$ commute.We will use the same inductive proof.
It is possible to show \cite{Ba} that elements $(X^N)_{k,l}$ are tensorial:
 $[X_{i,j},(X^N)_{k,l}]=\d_{k,j} (X^N)_{i,l} -\d_{i,l} (X^N)_{k,j}+\e_i \e_j(\d_{j,-l}(X^N)_{k,-i}-\d_{k,-i}(X^N)_{-j,l }).$

 Besides we will need the following proposition:

 {\large Proposition 3.}
$ (X^{M+1})_{i,j}=\sum\limits _{p=0}^{M+1} C_p \e_i \e_j  (X^p)_{-j,-i }$,

{\it where $C_p$ depend on the Casimir elements of the type 
$(X^k),k<p$, $C_{M+1}=(-1)^{M+1}$.}

{\it Proof .}
Proof of this proposition is inductive. Obviously,that it is true for $X^0_{i,j}=\d_{i,j}$. It is true by definition for $X^1_{i,j}=X_{i,j}$.
Let us suppose that it is true also for $M=M_0$. It is easy to show that

$(X^{M_0 + 1})_{i,j}= - \e_i \e_k (X^{M_0})_{k,j }X_{-k,-i}+\Tr (id ) (X^{M_0})_{i,j }- \e_i \e_j (X^{M_0})_{-j,-i }-\e_i \e_j (X^{M_0})(X^0)_{-j,-i }.$

So, from the assumption that the statement of proposition is true for $M=M_0$  it follows that it is true  also for $M_0 +1$. That proves the proposition.

 {\large Proposition 4.} $[(X^M)_{i,j}, (X^N)_{k,l}]=\sum\limits_{S=1}^{M}((X^{M+N-S})_ {i,l}  (X^{S-1}) _{k,j} -(X^{S-1})_{i,l} (X^{N+M-S})_{k,j})$

$-\sum\limits_{p=0}^N C_p (\sum\limits_{S=1}^{M}(\e_{-l}\e_k  (X^{M+p-S})_ {i,-k} ( X^{S-1}) _{-l,j} - \e_{-k}\e_l (X^{S-1})_{i,-k} (X^{p+M-S})_{-l,j}))$.

{\it Proof .}Using the Leibnitz rule we obtain that 

$[(X^M)_{i,j},(X^N)_{k,l}]=\sum\limits_{S=1}^M \sum\limits_{u,v \in I} (X^{S-1})_{i,u}[X_{u,v}, (X^N)_{k,l}] (X^{M-S})_{v,j}$.

Using tensorial nature of elements $(X^N)_{k,l}$ we have

$ [(X^M)_{i,j},(X^N)_{k,l}]
=\sum\limits_{S=1}^{M}((X^{M+N-S})_ {i,l} ( X^{S-1}) _{k,j} - (X^{S-1})_{i,l} (X^{N+M-S})_{k,j}$

$+ \sum\limits_{u \in I}(\e_u \e_{-l} (X^{S-1})_ {i,u} (X^{N})_{k,-u} (X^{M-S}) _{-l,j}  
- \e_u \e_{-k} (X^{S-1})_{i,-k} (X^{N})_{-u,l}  (X^{ M-S}_{u,j}))$.

Using Proposition 3 we obtain the statement of the Proposition 4.

We will prove the theorem 
 if we prove the equality
$[(AX^K), (AX^N)]=0,$ for all $N,K\in{\bf Z_+}$,
i.e. that $\sum\limits_{i,j,k,l\in I} A_{j,i}A_{l,k} [(X^K)_{i,j},(X^N)_{k,l}]=0 $.

We will need also equality
$\sum\limits_{i,j,k,l\in I} A_{j,k}A_{l,i} [(X^K)_{i,j}, (X^N)_{k,l}]=0 ,$
which will be interconnected with the first one in its inductive proof.
Let's proceed inductively, assuming, that these equalities  are valid for all $K<M$ and arbitrary $N\in{\bf  Z_+}$ . We have to draw out that these equalities are valid also for $K=M$.
These assumptions will be proved by 

{\large Proposition 5.}
{\it
 If $ A_{i,j}=  \pm \e_{j}\e_{-i} A_{-j, -i}$ then 
the next equalities hold:}

(i) $\sum\limits_{i,j,k,l\in I} A_{j,i}A_{l,k} [(X^M)_{i,j},(X^N)_{k,l}]=
\sum\limits_{S=1}^{M}\sum\limits_{i,j,k,l\in I} A_{j,i}A_{l,k}[(X^{S-1})_{i,j},(X^{N+M-S})_{k,l}]+$

$+\sum\limits_ {P=0}^{N} C_p (\sum\limits_ {S=1}^{M}\sum\limits_{i,j,k,l\in I} A_{j,k}A_{l,i}[(X^{S-1})_{i,j},(X^{P+M-S})_{k,l}])$.

(ii) $\sum\limits_{i,j,k,l\in I} A_{j,k}A_{l,i} [(X^M)_{i,j}, (X^N)_{k,l}]= \sum\limits_{S=1}^{M}\sum\limits_{i,j,k,l\in I} A_{j,i}A_{l,k}[(X^{S-1})_{i,j},(X^{N+M-S})_{k,l}]+$

$\pm \sum\limits_ {P=0}^{N}C_P \sum\limits_ {S=1}^{M}\sum\limits_ {Q=1}^{S}( C_Q \sum\limits_{i,j,k,l\in I} A_{l,i}A_{j, k}[(X^{Q-1})_{i,j}, (X^{P+M-S})_{k,l}])$.

{\it Proof:}
Let us prove part (i) of the Proposition.
Transforming its left-hand-side we obtain:

$\sum\limits_{i,j,k,l\in I} A_{j,i}A_{l,k} [ (X^M)_{i,j}, (X^N)_{k,l}]=
\sum\limits_{S=1}^{M}((AX^{S-1}AX^{M+N-S}) - (AX^{S-1}AX^{M+N-S}))-$

$-\sum\limits_{S=1}^{M}\sum\limits_{i,j,k,l\in I} A_{j,i}A_{l,k}[(X^{S-1})_{i,l}, (X^{N+M-S})_{k,j}]+$

$+\sum\limits_{p=0}^N C_p (\sum\limits_{S=1}^{M}(\e_{-l}\e_k  A_{j,i} A_{l,k}(X^{M+p-S})_ {i,-k} ( X^{S-1}) _{-l,j} - \e_{-k}\e_l A_{j,i} A_{l,k}(X^{p+M-S})_{-l,j} (X^{S-1})_{i,-k}))$

$-\sum\limits_ {P=0}^{N} C_p (\sum\limits_ {S=1}^{M}\sum\limits_{i,j,k,l\in I}\e_{-l}\e_k A_{j,i}A_{l,k}[(X^{S-1})_{i,-k}, (X^{P+M-S})_{-l,j}])$.

First summand is identically zero. It is easy to see that after relabelling of the indices  in the second summand $j\leftrightarrow l$ we will have it in the form of (i). If we relabell indices $i\leftrightarrow -l$, $j\leftrightarrow -k$ in the second part of the third summand we will have it in the following form

$ \sum\limits_{p=0}^N C_p \sum\limits_{S=1}^{M}(\e_{-l}\e_k  A_{j,i} A_{l,k}(X^{M+p-S})_ {i,-k} ( X^{S-1}) _{-l,j} - \e_{j}\e_{-i} A_{-i, -j} A_{-k,-l}(X^{p+M-S})_{i,-k}( X^{S-1})_{-l,j})$

from which it is easy to see that this summand is identically zero if and only if
$A_{i,j}=  \pm \e_{j}\e_{-i} A_{-j, -i}.$
Under this condition after relabelling the indices $j\rightarrow l$,$-l\rightarrow k$,
$k\rightarrow -j$ it is easy to see that fourth summand coinsides with the second summand in the right-hand-side of (i). We have proved (i).
 By analogous direct calculations, using Proposition 3, one can verify validity of (ii).
Now to prove the theorem it remains only to show  that

$\sum\limits_{i,j,k,l\in I} A_{j,i}A_{l,k} [X _{i,j}, (X^N)_{k,l}]=0$, and
$\sum\limits_{i,j,k,l\in I} A_{j,k}A_{l,i} [X _{i,j}, (X^N)_{k,l}]=0$.

First equality follows from the fact that $(AX) \in g_A$.
Let us prove the second equality. Transforming its left-hand-side we obtain

 $\sum\limits_{i,j,k,l\in I} A_{j,k}A_{l,i} [X _{i,j},(X^N)_{k,l}]=
\sum\limits_{i,j,k,l\in I} A_{j,k}A_{l,i}((\d_{k,j} (X^N)_{i,l} -\d_{i,l}( X^N)_{k,j})+$

$+\e_i \e_j(\d_{j,-l}(X^N)_{k,-i}-\d_{k,-i}(X^N)_{-j,l }))= \sum\limits_{i,j,k,l\in I}( \e_i\e_{-l}A_{l,i}A_{- l,k}(X^N)_{k,-i}-$

$- \e_i \e_{j}A_{j,-i} A_{l,i}(X^N)_{-j,l})
=(A^2 X^N) - (A^2 X^N)=0$.

That  completes the proof of the theorem.

\section{ Proof of  Theorem 3.}
Taking into account commutativity of the algebras $F_A$ and $C_{g_A}$ provided by  Theorems 1-2 and the assumption of Theorem 3 we can use 
the result of Vinberg  \cite {V} that the degree of transcendency of any commutative subalgebra in an associative filtered Lie algebra is equal to the degree of transcendency of correspondent subalgebra in the associated graded algebra.
That is why we may consider commutative (with respect to the Poisson
brackets ) subalgebra in $P(g^*)$ \footnote { One can identify $g$ and $g^*$ in the semisimple case.} instead any commutative subalgebra in $U(g)$.
Let
$ S^M (X)=(X^M)=\sum\limits_ { i_1,...i_{M} \in I} X_{i_1,i_2}X_{i_3,i_4}...X_{i_{M}, i_1 }$
be a Casimir function in $P(g^*)$.
We still consider all classical simple algebras as a subalgebras in $gl(n)$ and denote coordinate functions  on $g^*$ by the same letters as basic elements of $g$.
It is easy to see that  commutative subagebra $F_A (g)$ coincides with the subalgebra of the integrals of the degree not higher than one in the parameters of the shift of the Mishchenko-Fomenko algebra.
Let us explain its construction. Let $A \in g^*$ be any constant covector. Then one can define functions $S_A^{k,M}(X), k<M$ from the expansion:
$S^M (X+\lambda A) =\sum\limits_{k=1}^M  S_A^{k,M} (X) \lambda ^k$.
Let $F(A,X)$ be commutative algebra generated by all functions $S_A^{k,M}(X), k<M, M \in Z_+$. Following \cite{V} we call it
Mishchenko-Fomenko algebra.
We will need the following Lemma:

 {\large Lemma 2.}
{\it Let $O_{sing}$ be the singular semisimple coadjoint orbit in $g^*$ consisting of the elements of $g^*$ of the second matrix rank.Then for independent generators of  the Mishchenko-Fomenko algebra  restricted to this orbit could be chosen those ones  having degrees not higher than 2 in the coordinates of $g^*$.}

{\it Sketch of the  proof  \footnote{For more detailed proof look \cite{BS}.}.} 
Let $\lbrace P^M (X) \rbrace$ be a complete set of invariant functions chosen in the form of coeficients of the characteristic polinomial. Let $P_A^{k,M}(X), k<M$ be the set of "integrals" obtained by the described above procedure of the "shift of the argument".
It is easy to show that $P_A^{k,M}(X)$ are linear combinations of the minor determinants of the matrix $X$ of the order $M-k$. Hence after restriction to the $O_{sing}$ only those with $M-k<3$ (i. e. those of the degree not higher than 2 in the $X$ variables) are not zero.
Statement of the proposition now follows from the fact that $\lbrace P^M (X) \rbrace$ and $\lbrace S^M (X) \rbrace$ and, hence, $\lbrace P_A^{k,M}(X), k<M \rbrace$
and $\lbrace S_A^{k,M}(X), k<M \rbrace$ are connected by a bipolinomial correspondence.
\vskip 2mm
We can consider $S_A^{k,M}(X) $ as functions depending on the $X$ and $A$. 
Let $d_X (.)$ denotes the "matrix" gradients of  functions $S_A^{k,M}(X) $ with the values in $g$. Let
 $d_A (.)$ denotes the matrix gradients of  functions $S_A^{k,M}(X) $
considered as functions of  parameters $A$. Functional independence of functions in $ P(g^*)$ (as functions of  $X$) is equivalent to the  linear independence of their gradients.
Let $D_X(f)$( $D_A (f)$) denotes the space of matrix gradients over variables $X$ (parameters $A$) of all generators of some subalgebra $f\subset F(A,X)$.
There exists a kind of duality between integrals from the Mishchenko-Fomenko algebra as functions of $X$ and $A$ variables \cite {Br}:
$d_X (S_A^{k,M}(X)) = d_A (S_X^{M-k-1,M}(A)), $ 
\hskip 2 mm
$D_X (F(A,X)) = D_A (F (X,A)).$\footnote{ Note that in the right-hand-side of these relations $X_{i,j}$ play the role of the parameters of the shift and $A_{i,j}$ play the role of the coordinates of the algebra.}
From the first equality follows that
$  D_X (F_A(X))= D_A (F^2_X(A))$, where $F^2_X(A)$ is a subalgebra in $F(X, A)$ of the first and the second order in parameters $A$. 
From the Lemma 2 follows that $D_A (F^2_X (A))= D_A (F(X,A))$.
Hence we obtain
$  D_X (F_A(X))= D_A (F (X,A))= D_X (F(A,X))$.
As it was shown in \cite {Br}:
$\dim (D_X (F(A,X))\bigcap [A,g])= 1/2 \dim [A,g]$.
Hence we obtain:
 $\dim (D_X (F_A(X))\bigcap [A,g])= 1/2  \dim [A,g]$.
Taking into account, that
$[A,g] \perp g_A$,  $ \dim[A,g]$=$ \dim g-\dim g_A$, and
$\ind g_A=\ind g=\rank g$ 
for semisimple matrix $A$, we obtain the statement of the theorem.

\end{document}